\documentclass[12pt,reqno]{amsart}
\usepackage{amssymb,latexsym,graphicx,amscd,upgreek}
\usepackage{lscape}
\usepackage[all]{xy}
\usepackage{color}
\usepackage{mathrsfs}
\usepackage{setspace}
\usepackage{cases}
\usepackage{bbm}
\usepackage{hyperref}
\hypersetup{
    colorlinks,
    linkcolor={blue},
    citecolor={blue},
    urlcolor={blue!80!black}
}

\newtheorem{Thm}{Theorem}[section]
\newtheorem{Lem}[Thm]{Lemma}
\newtheorem{Cor}[Thm]{Corollary}
\newtheorem{Prop}[Thm]{Proposition}

\theoremstyle{definition}
 \newtheorem{Rem}[Thm]{Remark}
 \newtheorem{Def}[Thm]{Definition}

\newcommand{\Example}{\noindent{\bf Example}}
\newcommand{\Examples}{\noindent{\bf Examples}}

\newcommand{\Z}{\mathbb{Z}}

\newcommand{\df}{\colon}

\newcommand{\cL}{{\mathcal L}}

\newcommand{\scL}{{\mathscr L}}
\newcommand{\scR}{{\mathscr R}}

\newcommand{\cP}{{\mathcal P}}
\newcommand{\cR}{{\mathcal R}}

\newcommand{\m}{\mathfrak{m}}

\newcommand{\vep}{\varepsilon}

\newcommand{\cpQ}{\langle\!\langle Q \rangle\!\rangle}

\newcommand{\rk}{\operatorname{rk}}

\newcommand{\md}{\operatorname{mod}}

\newcommand{\Jacobi}{\cP}

\newcommand{\gldim}{\operatorname{gl.dim}}
\newcommand{\pdim}{\operatorname{proj.dim}}
\newcommand{\idim}{\operatorname{inj.dim}}

\newcommand{\dimv}{\underline{\dim}}

\newcommand{\soc}{\operatorname{soc}}
\newcommand{\rad}{\operatorname{rad}}
\newcommand{\tp}{\operatorname{top}}

\newcommand{\cogen}{\operatorname{Cogen}}

\newcommand{\GP}{\operatorname{GP}}

\newcommand{\Hom}{\operatorname{Hom}}
\newcommand{\Ext}{\operatorname{Ext}}

\newcommand{\ext}{\operatorname{ext}}
\newcommand{\End}{\operatorname{End}}

\newcommand{\Aut}{\operatorname{Aut}}

\newcommand{\Ima}{\operatorname{Im}}
\newcommand{\Ker}{\operatorname{Ker}}
\newcommand{\Coker}{\operatorname{Cok}}

\newcommand{\ov}{\overline}

\newcommand{\bil}[1]{\langle #1\rangle}

\newcommand{\bsm}{\begin{smallmatrix}}
\newcommand{\esm}{\end{smallmatrix}}

\newcommand{\bbsm}{\left[\begin{smallmatrix}}
\newcommand{\besm}{\end{smallmatrix}\right]}

\newcommand{\bbm}{\begin{matrix}}
\newcommand{\ebm}{\end{matrix}}

\newcommand{\bLam}{{\overline{\Lambda}}}

\begin{document}
%%%%%%%%%%%%%%%%
%%%%%%%%%%%%%%%%

%\today
\date{10.11.2026}

%%%%%%%%%%%%%%%%%%%%%%%%%%%%%%%%%%%%%%%%%%%%%%
\title{$E$-symmetric algebras II}
%%%%%%%%%%%%%%%%%%%%%%%%%%%%%%%%%%%%%%%%%%%%%%%

\author{Grzegorz Bobi\'nski}
\address{Grzegorz Bobi\'nski\newline
Faculty of Mathematics and Computer Science\newline
Nicolaus Copernicus University\newline
ul. Chopina 12/18\newline
87-100 Toru\'n\newline
Poland}
\email{gregbob@mat.umk.pl}

%\author{Christof Gei{\ss}}
%\address{Christof Gei{\ss}\newline
%Instituto de Matem\'aticas\newline
%Universidad Nacional Aut{\'o}noma de M{\'e}xico\newline
%Ciudad Universitaria\newline
%04510 M{\'e}xico D.F.\newline
%M{\'e}xico}
%\email{christof.geiss@im.unam.mx}

\author{Carolin Hartung}
\address{Carolin Hartung\newline
Mathematisches Institut\newline
Universit\"at Bonn\newline
Endenicher Allee 60\newline
53115 Bonn\newline
Germany}
\email{hartung@math.uni-bonn.de}

%\author{Daniel Labardini-Fragoso}
%\address{Daniel Labardini-Fragoso\newline
%Dipartimento di Matematica 'Tullio Levi-Civita'\newline
%Università degli Studi di Padova\newline
%Via Trieste, 63 - 35121 Padova\newline
%Italy
%}
%\email{daniel.labardinifragoso@unipd.it}

\author{Jan Schr\"oer}
\address{Jan Schr\"oer\newline
Mathematisches Institut\newline
Universit\"at Bonn\newline
Endenicher Allee 60\newline
53115 Bonn\newline
Germany}
\email{schroer@math.uni-bonn.de}


\begin{abstract}
We introduce generalized Jacobian algebras and show that a 
finite-dimensional algebra is $E$-symmetric if and only if it is
generalized Jacobian.
We also characterize $E$-symmetric algebras in terms of its simple modules.
We show that 
Enomoto's counterexample to the Periodicity Conjecture 
is a generalized Jacobian algebra, which is not Jacobian in the classical sense.
\end{abstract}

%%%%%%%%%%%%%%%%
\maketitle
%%%%%%%%%%%%%%%%

\setcounter{tocdepth}{1}
\numberwithin{equation}{section}
\tableofcontents

\parskip2mm
%\parindent0pt

%%%%%%%%%%%%%%%%%%%%%%%%%%%%%%%%%%%%%%%
%%%%%%%%%%%%%%%%%%%%%%%%%%%%%%%%%%%%%%%

\section{Introduction and main results}\label{sec:1}

%%%%%%%%%%%%%%%%%%%%%%%%%%%%%%%%%%%%%%%
%%%%%%%%%%%%%%%%%%%%%%%%%%%%%%%%%%%%%%%

Throughout, $K$ denotes a field, and let
$A = KQ/I$ be a finite-dimensional $K$-algebra, where $KQ$ is the path algebra of a quiver $Q$ and 
$I$ is an admissible ideal of $KQ$.
Let $\md(A)$ be the category of finite-dimensional left $A$-modules.

%%%%%%%%%%%%%%%%%%%%%%%%%%%%%
\subsection{$E$-invariants}
%%%%%%%%%%%%%%%%%%%%%%%%%%%%%

For $M,N \in \md(A)$, we define
$$
\hom_A(M,N) := \dim \Hom_A(M,N).
$$

Let
\begin{align*}
E(M,N) &:= \hom_A(N,\tau_A(M)),
\\
E^-(M,N) &:= \hom_A(\tau_A^{-}(N),M),
\end{align*}
where $\tau_A$ denotes the Auslander-Reiten translation for $A$.
Recall that the Auslander-Reiten formula implies that there are functorial
isomorphism
$$
\ov{\Hom}_A(N,\tau_A(M)) = \underline{\Hom}_A(\tau_A^{-}(N),M),
$$
where $\ov{\Hom}$ and $\underline{\Hom}$ denote stable homomorphism spaces, see for example \cite[Section~IV.4]{ARS97} for details.

By \cite[Proposition~1.5]{BS26}, we have
$E(M,N) = E^-(M,N)$ for all $M,N \in \md(A)$ if and only if $A$ is hereditary.

Let
$$
E(M) := E(M,M) \text{\quad and \quad}
E^-(M) := E^-(M,M).
$$
These numbers are called \emph{$E$-invariants} of $M$ and
play an important role in Derksen--Weyman--Zelevinsky's
\cite{DWZ08, DWZ10} representation theoretic approach to 
Fo\-min--Zelevinsky cluster algebras.

We call $A$ an \emph{$E$-symmetric algebra}
if $$
E(M) = E^-(M)
$$ 
for all $M \in \md(A)$. 
The study of $E$-symmetric algebras was initiated in \cite{BS26}.

%%%%%%%%%%%%%%%%%%%%%%%
\subsection{Generalized Jacobian algebras}\label{subsec:1.2}
%%%%%%%%%%%%%%%%%%%%%%%)
Let $A = KQ/I$ as before. 
Let $\m$ be the arrow ideal of $KQ$.
We know that there is some $N_I \ge 2$ such that
$$
\m^{N_I} \subseteq I \subseteq \m^2.
$$

For arrows $a,b \in Q_1$ with $t(a) = k = s(b)$, we write their composition
as $ba$.
$$
\xymatrix{
s(a) \ar[r]^a  & k \ar[r]^b & t(b)
}
$$

For $i \in Q_0$, let $P(i) := Ae_i$ (resp. $I(i) := D(e_iA)$) be the indecomposable projective (resp. injective) $A$-module associated to $i$, and let $S(i) \cong \tp(P(i)) \cong \soc(I(i))$ be the corresponding simple
$A$-module.
Here $e_i$ denotes the path of length $0$ associated to $i$.

For each $a \in Q_1$, we define $K$-linear \emph{arrow removal maps} 
$$
r_a, l_a\df KQ \to KQ
$$ 
which are defined on
paths $p$ by
$$
r_a(p) :=
 \begin{cases}
 q&:\;  \text{if $p = qa$ for some path $q$},
\\
0 &:\;  \text{otherwise},
\end{cases}
$$
and
$$
l_a(p) :=
 \begin{cases}
 q&:\; \text{if $p=aq$ for some path $q$},
\\
0&:\;  \text{otherwise}.
 \end{cases}
$$
In particular,
\[
r_a(a) = e_{t(a)},\qquad 
l_a(a)=e_{s(a)},
\text{\quad and \quad} l_a(e_i) = r_a(e_i) = 0
\]
for all $i \in Q_0$.

\Example:
For a path $p = abc$ with $a,b,c \in Q_1$, we have $r_c(p) = ab$ and $l_a(p) = bc$.

Note that
for each $x \in e_j \m e_i$, we have  
$$
x = \sum_{a \in Q_1 : \; t(a) = j}  al_a(x)
= \sum_{b \in Q_1 : \; s(b) = i}  r_b(x)b.
$$

For $x \in KQ$, define $\ov{x} := x+I \in A$.

We get $K$-linear maps
$$
\ov{l_a}, \ov{r_a}\df KQ \to A
$$
defined by $x \mapsto \overline{l_a(x)}$ and
$x \mapsto \overline{r_a(x)}$, respectively.

\begin{Def}\label{def:genhyper}
A \emph{generalized hyperpotential} for $A = KQ/I$ consists of
an integer $N_I \geq 2$ with $\m^{N_I} \subseteq I$ and 
a set
\[
W = \{ \rho_a \in e_{s(a)}Ie_{t(a)} \mid a\in Q_1 \}
\]
such that the following hold:

\begin{itemize}\itemsep2mm

\item[(h1)] 
We have
\begin{equation}\label{eq:G1}
I = I(W) + \m^{N_I+1},
\end{equation}
where $I(W)$ is the ideal of $KQ$ generated by $W$.

\item[(h2)]
For each $k \in Q_0$, let 
\[
P_{k,{\rm in}} := \bigoplus_{a \in Q_1 : \; t(a)=k} P(s(a))
\text{\quad and \quad}
P_{k,{\rm out}} := \bigoplus_{b \in Q_1 : \; s(b)=k} P(t(b)).
\]
Define $A$-module homomorphisms
$$
\scL_k  = (\scL_{k,b,a})_{b,a},\;
\scR_k    = (\scR_{k,b,a})_{b,a}\df P_{k,{\rm in}}\to P_{k,{\rm out}}
$$ 
by
\begin{align*}
\scL_{k,b,a}\df P(s(a)) &\to P(t(b)) 
&
\scR_{k,b,a}\df P(s(a)) &\to P(t(b)) 
\\
x_a &\mapsto x_a \overline{l_a}(\rho_b)
&
x_a &\mapsto x_a \overline{r_b}(\rho_a).
\end{align*}
We think of $\scL_k$ and $\scR_k$ as matrices with entries
$\scL_{k,b,a}$ and $\scR_{k,b,a}$, where $b$ is a row index
and $a$ is a column index.
%We therefore have
%\[
% \begin{aligned}
%\scL_k((x_a)_a)
% &= \left(
%     \sum_{a \in Q_1 : \; t(a)=k}x_a\,\overline{l_a}(\rho_b)
%   \right)_{b \in Q_1 : \; s(b)=k},
%\\
%\scR_k((x_a)_a)
% &= \left(
%     \sum_{a \in Q_1 : \; t(a)=k}x_a\,\overline{r_b}(\rho_a)
%   \right)_{b \in Q_1 : \; s(b)=k}.
% \end{aligned}
%\]
We demand that for each $k \in Q_0$, there exist
\[
h_{k,{\rm in}} \in \Aut_A(P_{k,{\rm in}})
\text{\quad and \quad}
h_{k,{\rm out}} \in \Aut_A(P_{k,{\rm out}})
\]
such that
\begin{equation}\label{eq:G2}
\scL_k = h_{k,{\rm out}}^{-1} \circ \scR_k \circ h_{k,{\rm in}}.
\end{equation}
\end{itemize}
\end{Def}

\begin{Def}\label{def:genjacobian}
We call $A$ a \emph{generalized Jacobian algebra} if there exists
a generalized hyperpotential for $A$. 
\end{Def}

\begin{Rem}
In condition (h1), we work with $\m^{N_I+1}$ instead of $\m^{N_I}$, because 
$\overline{l_a}(p) = \overline{r_a}(p) = 0$ for all $p \in \m^{N_I+1}$.
In other words, the maps $\ov{l_a}$ and $\ov{r_a}$ are not affected by the extra relations coming from $\m^{N_I+1}$.
\end{Rem}

\begin{Rem}
Let $a,b \in Q_1$.
We have
$$
l_a(\rho_b),\; r_b(\rho_a)
\in
e_{s(a)}KQe_{t(b)}
\text{\quad and \quad}
\overline{l_a}(\rho_b),\; 
\overline{r_b}(\rho_a)
\in
e_{s(a)} A e_{t(b)}. 
$$
Note that for
$t(a) \not= s(b)$, we have
$r_b(\rho_a) = l_a(\rho_b) = 0$.
For
$t(a) = k = s(b)$,
the following picture illustrates the situation:
$$
\xymatrix@+2ex{
s(a) \ar[r]^a 
& k \ar[r]^b \ar@{-->}@/^4ex/[l]^{\rho_a}
& t(b) \ar@{-->}@/^4ex/[l]^{\rho_b}
&
s(a) \ar[r]^a & k \ar[r]^b
& t(b) \ar@{-->}@/^4ex/[ll]^{r_b(\rho_a)}
\ar@{-->}@/_4ex/[ll]_{l_a(\rho_b)}
}
$$
The main idea of the definition of a generalized hyperpotential
is that the two $A$-module homomorphisms
$$
\scL_k,\;
\scR_k\df P_{k,{\rm in}}\to P_{k,{\rm out}}
$$
satisfy the compatibility condition 
$$
\scL_k = h_{k,{\rm out}}^{-1} \circ \scR_k \circ h_{k,{\rm in}}.
$$
\end{Rem}

\begin{Rem}
A generalized hyperpotential as defined above is a
\emph{hyperpotential} as defined by Ladkani \cite{L14}
if and only if
$$
\scL_k = \scR_k
$$
for all $k \in Q_0$.
Hyperpotentials generalize the notion of \emph{potentials} as defined by Derksen--Weyman--Zelevinsky \cite{DWZ08,DWZ10}.
\end{Rem}

For $i \ge 0$, let $KQ_i$ be the subspace of $KQ$ with basis
the set $Q_i$ of paths of length $i$ in $Q$.
For $m \ge 0$, let
$$
KQ_{\le m} := \bigoplus_{0 \le i \le m} KQ_i.
$$

\begin{Rem}
In the definition of a generalized hyperpotential, one can assume without loss of generality that $$\rho_a \in KQ_{\le N_I}.$$
To be more precise, one can 
write $\rho_a = \rho_a' + \rho_a''$ with
$\rho_a' \in KQ_{\le N_I}$ and $\rho_a'' \in \m^{N_I+1}$
and then replace $\rho_a$ by $\rho_a'$. 
Condition (h1) is then still satisfied, and the maps $\scL_k$ and 
$\scR_k$ do not change.
\end{Rem}

\begin{Rem}
It follows from a general construction of the first three terms of
projective resolutions of simple modules as explained in
\cite[Section~5.4]{GHLS26}, 
that
$\scR_k$ is a 
(not necessarily minimal) 
projective presentation
of $\rad(P(k))$ for each $k \in Q_0$.
The compatibility condition in condition (h2) implies then
$$
\Coker(\scL_k) \cong \Coker(\scR_k) \cong \rad(P(k)).
$$
\end{Rem}

\begin{Rem}
As the name suggests, the class of
generalized Jacobian algebras contains the class of Jacobian algebras defined by potentials or hyperpotentials.
By \cite{DWZ10} (for potentials) and \cite{GHLS26} (for hyperpotentials), Jacobian algebras $\Jacobi_K(Q,W)$ are $E$-symmetric.
\end{Rem}

%%%%%%%%%%%%%%%%%%%
\subsection{Main result}
%%%%%%%%%%%%%%%%%%%
Recall that $A$ is a \emph{$1$-Iwanaga-Gorenstein algebra}, 
if $$
\pdim(D(A_A)) \le 1
\text{\quad and \quad}
\idim({_A}A) \le 1.
$$
In this case, the
category $\GP(A)$ of Gorenstein-projective $A$-modules consists of
all $A$-modules which are isomorphic to a submodule of some projective $A$-module.

\begin{Thm}\label{thm:main1}
The following are equivalent:
\begin{itemize}\itemsep2mm

\item[(i)]
$A$ is $E$-symmetric;

\item[(ii)]
For each simple $A$-module $S$, we have 
an isomorphism
$$
\tau_A(\Omega_A(S)) \cong \Omega_A^{-1}(S)
$$ 
in the stable category $\overline{\md}(A)$;

\item[(iii)]
$A$ is a $1$-Iwanaga-Gorenstein algebra, and for each simple $A$-module $S$ and $M \in \GP(A)$, we have an isomorphism
$$
\Ext_A^1(M,S) \cong D\!\Ext_A^2(S,M)
$$
which is functorial in $M$;

\item[(iv)]
$A$ is a generalized Jacobian algebra.

\end{itemize}
\end{Thm}

\begin{Rem}
Condition (ii) in Theorem~\ref{thm:main1} can be reformulated as follows:
\begin{itemize}\itemsep2mm

\item[]
\emph{
For each simple $A$-module $S$, we have 
an isomorphism
$$
\tau_A(\Omega_A(S)) \cong \Omega_A^{-1}(S)/C_S
$$ 
in $\md(A)$, where $C_S$ is a maximal injective direct summand
of $\Omega_A^{-1}(S)$.}

\end{itemize}
\end{Rem}

\begin{Rem}
Let $A = KQ/I$ be a monomial algebra.
Then it follows from
\cite[Theorem~1.3]{BS26} and \cite[Theorem~1]{GE21}
that
the following are equivalent:
\begin{itemize}\itemsep2mm

\item[(i)]
$A$ is $E$-symmetric;

\item[(ii)]
$A$ is $1$-Iwanaga-Gorenstein, and for each $M \in \GP(A)$ and $N \in \md(A)$, we have an isomorphism
$$
\Ext_A^1(M,N) \cong D\!\Ext_A^2(N,M)
$$
which is functorial in $M$ and $N$;

\item[(iii)]
$A$ is a Jacobian algebra.

\end{itemize}
\end{Rem}

\begin{Rem}
Theorem~\ref{thm:main1} can be generalized from 
finite-dimensional algebras 
$A = KQ/I$ to
arbitrary basic pseudocompact algebras
$\bLam = K\cpQ/\overline{I}$. 
The proof uses the same arguments, but the whole setup is more technical. 
We will do this in a separate publication \cite{BHS26}. 
\end{Rem}

We discuss now two easy examples. A more elaborate example
can be found in Section~\ref{sec:6}.

\Examples:
\begin{itemize}\itemsep2mm

\item[(i)]
Let $A = KQ/I$, where $Q$ is the quiver
\[
\xymatrix@C=3mm@R=6mm{
  & 2 \ar[dr]^b & \\
  1 \ar[ur]^a & & 3 \ar[ll]^c && 4 \ar[ll]^d
}
\]
and $I$ is generated by $\{  cb,\; ac,\;  ba \}$.
We display the relevant $A$-modules:
$$
\begin{tabular} {c|c|c|c}
$S(i)$ & $\Omega_A(S(i))$ & $\tau_A(\Omega_A(S(i)))$ 
& $\Omega_A^{-1}(S(i))$
\\\hline
1 & $\begin{matrix} 2 \end{matrix}$ &$\begin{matrix} 4\\3 \end{matrix}$&
$\begin{matrix} 4\\3 \end{matrix}$
\\\hline
2 &$\begin{matrix} 3 \end{matrix}$&$\begin{matrix} 1 \end{matrix}$&
$\begin{matrix} 1 \end{matrix}$
\\\hline
3 &$\begin{matrix} 1 \end{matrix}$&$\begin{matrix} 2 \end{matrix}$&
$\begin{matrix} 2 \oplus 4 \end{matrix}$
\\\hline
4 &$\begin{matrix} 3\\1 \end{matrix}$&$\begin{matrix} 0 \end{matrix}$&
$\begin{matrix} 0 \end{matrix}$
\end{tabular}
$$
Each number $i$ stands here for a composition factor $S(i)$.
Note that $4 = S(4) \cong I(4)$.
Using the criterion in Theorem~\ref{thm:main1}, 
it follows that
$A$ is $E$-symmetric.
Let $N_I := 3$.
Then 
$$
W := \{ \rho_a = cb,\; \rho_b = ac,\; \rho_c =  ba,\; \rho_d = 0\}
$$ 
(together with $N_I$)
is a generalized hyperpotential for $A$.
For example, for $k=3$, the maps $\cL_k$ and $\cR_k$ look as follows:
\begin{align*}
\cL_3  = (\ov{l_b}(\rho_c),\ov{l_d}(\rho_c)) = (a,0) \df P(2) \oplus P(4) &\to P(1),
\\
\cR_3  = (\ov{r_c}(\rho_b),\ov{r_c}(\rho_d)) = (a,0) \df P(2) \oplus P(4) &\to P(1),
\end{align*}
where the matrix entries mean right multiplication by (the residue classes modulo $I$ of) the displayed elements.
For all $1 \le k \le 4$, we have $\cL_k = \cR_k$.

\item[(ii)]
Let $A = KQ/I$, where $Q$ is the quiver
\[
\xymatrix{
1  & 2 \ar[l]_a & 3 \ar[l]_b 
}
\]
and $I$ is generated by $ab$. Here are
the relevant $A$-modules:
$$
\begin{tabular} {c|c|c|c}
$S(i)$ & $\Omega_A(S(i))$ & $\tau_A(\Omega_A(S(i)))$ 
& $\Omega_A^{-1}(S(i))$
\\\hline
1 & $\begin{matrix} 0 \end{matrix}$ &$\begin{matrix} 0 \end{matrix}$&
$\begin{matrix} 2 \end{matrix}$
\\\hline
2 &$\begin{matrix} 1 \end{matrix}$&$\begin{matrix} 0 \end{matrix}$&
$\begin{matrix} 3 \end{matrix}$
\\\hline
3 &$\begin{matrix} 2 \end{matrix}$&$\begin{matrix} 1 \end{matrix}$&
$\begin{matrix} 0 \end{matrix}$
\end{tabular}
$$
Thus $A$ is not $E$-symmetric.

\end{itemize}

For the class of selfinjective algebras, Theorem~\ref{thm:main1}
yields the following characterization of $E$-symmetry.

\begin{Cor}\label{cor:main2}
Let $A$ be selfinjective.
Then the following are equivalent:
\begin{itemize}\itemsep2mm

\item[(i)]
$A$ is $E$-symmetric;

\item[(ii)]
For each non-projective simple $A$-module $S$, we have
$$
\Omega_A^4(S) \cong \nu_A^{-1}(S).
$$

\end{itemize}
\end{Cor}

Here 
$\Omega_A$ and $\nu_A$
are the syzygy functor and the Nakayama functor, respectively.

We also obtain a characterization of local $E$-symmetric algebras:

\begin{Cor}\label{cor:main4}
Let $A$ be local.
Then the following are equivalent:
\begin{itemize}\itemsep2mm

\item[(i)]
$A$ is $E$-symmetric;

\item[(ii)]
$A$ is selfinjective, and if the simple $A$-module $S$ is
non-projective, then
$$
\Omega_A^4(S) \cong S.
$$

\end{itemize}
\end{Cor}

For $0 \not= q \in K$, let $A_q$ be the generalized Jacobian algebra
which arises as a deformation of the Jacobian algebra associated to the Markov quiver with the Labardini-Fragoso potential.
For a precise definition, we refer to Section~\ref{sec:6}.
For $q \not= 1$, the algebra $A_q$ is a counterexample to 
the Periodicity Conjecture, as reported in \cite{E26}.

\begin{Thm}\label{thm:main3}
For $q \not= 1$,
the algebra $A_q$ is $E$-symmetric, but it is not isomorphic to a Jacobian algebras defined by
some hyperpotential.
\end{Thm}

%%%%%%%%%%%%%%%
\subsection{Conventions}
%%%%%%%%%%%%%%%
For maps $f\df U \to V$ and $g\df V \to W$, we denote their
composition by $gf\df U \to W$.

For a linear map $f$ of vector spaces, let $\rk(f) := \dim \Ima(f)$
be its \emph{rank}.

%%%%%%%%%%%%%%%%%%%%%%%%%%%%%
\subsection{Organization of the article}
%%%%%%%%%%%%%%%%%%%%%%%%%%%%%

In
Section~\ref{sec:1}, we define generalized Jacobian algebras
and state the main result 
Theorem~\ref{thm:main1} (= Theorems~\ref{thm:mainpart1} and \ref{thm:mainpart2}).
Section~\ref{sec:IG} is a short recollection of known results on
$1$-Iwanaga-Gorenstein algebras.
In
Section~\ref{sec:2}, we
recall results on standard complexes (in the sense of
\cite{GHLS26}) and $g$-vectors.
It also contains a characterization of $E$-symmetric algebras in terms of the ranks of the maps  $M(\gamma_{k,\cR}')$ and
$M(\gamma_{k,\cR}'')$ from the standard complexes, and in terms of $g$-vectors, see Proposition~\ref{prop:EsymmRank}.
The main result of
Section~\ref{sec:3} is 
Theorem~\ref{thm:mainpart1}, which forms one part of 
Theorem~\ref{thm:main1}.
In
Section~\ref{sec:4}, we 
show that the ranks of $M(\gamma_{k,\cR}')$ and
$M(\gamma_{k,\cR}'')$ 
equal the ranks of $\Hom_A(\scL_k,M)$ and $\Hom_A(\scR_k,M)$, respectively, see Corollary~\ref{cor:h2}.
This allows us to finalize the proof of
Theorem~\ref{thm:main1}.
In
Section~\ref{sec:5}, we prove
Corollary~\ref{cor:main2} (= Corollary~\ref{cor:main2b})
and Corollary~\ref{cor:main4} (= Corollary~\ref{cor:main4b}),
which characterize $E$-symmetric selfinjective algebras
and $E$-symmetric local algebras, respectively.
Finally,
Section~\ref{sec:6} contains the proof of Theorem~\ref{thm:main3}.

%%%%%%%%%%%%%%%%%%%%%%%%%%%%%%%%
%%%%%%%%%%%%%%%%%%%%%%%%%%%%%%%%

\section{Iwanaga-Gorenstein algebras and Gorenstein-projective modules}
\label{sec:IG}

%%%%%%%%%%%%%%%%%%%%%%%%%%%%%%%%
%%%%%%%%%%%%%%%%%%%%%%%%%%%%%%%%

%%%%%%%%%%%%%%%%%%%%%%%%%%%%%%%%
\subsection{Iwanaga-Gorenstein algebras}
%%%%%%%%%%%%%%%%%%%%%%%%%%%%%%%%
Recall that
an algebra $A$ is a \emph{$1$-Iwanaga-Gorenstein algebra} if
$$
\pdim(D(A_A)) \le 1 
\text{\quad and \quad}
\idim({_A}A) \le 1.
$$
In this case, for $M \in \md(A)$, the following are equivalent:
\begin{itemize}\itemsep2mm

\item[(i)]
$\pdim(M) \le 1$;

\item[(ii)]
$\pdim(M) < \infty$;

\item[(iii)]
$\idim(M) \le 1$;

\item[(iv)]
$\idim(M) < \infty$.

\end{itemize}

%%%%%%%%%%%%%%%%%%%%%%%%%%%%%%%%%
\subsection{Gorenstein-projective modules}
%%%%%%%%%%%%%%%%%%%%%%%%%%%%%%%%%
For $M \in \md(A)$, let
$$
\cogen(M) := \{ U \in \md(A) \mid \text{there is a monomorphism $U \to M^t$ for some $t \ge 1$} \}.
$$
For a $1$-Iwanaga-Gorenstein algebra $A$,   
\begin{align*}
\GP(A) &:= \{ M \in \md(A) \mid \Ext_A^1(M,{_A}A) = 0 \}
\\
&= \cogen({_A}A)
\end{align*}
is the subcategory of \emph{Gorenstein-projective} $A$-modules.
In particular, for all $N \in \md(A)$, we have $\Omega_A(N) \in \GP(A)$.
The category $\GP(A)$
is a Frobenius subcategory of $\md(A)$ whose projective-injective objects
are the projective $A$-modules.

%%%%%%%%%%%%%%%%%%%%%%%%%%%%%%%%%%%%%%%%%%%%%
\subsection{Calabi-Yau condition on Gorenstein-projectives}
%%%%%%%%%%%%%%%%%%%%%%%%%%%%%%%%%%%%%%%%%%%%%
Let $A$ be a $1$-Iwanaga-Gorenstein algebra.
Then $A$
is \emph{stably $3$-Calabi-Yau} if there are functorial isomorphisms
$$
\Ext_A^1(M,N) \cong D\!\Ext_A^2(N,M) 
$$
for all $M,N \in \GP(A)$.
(This is a slightly weaker version of the definition used for example in \cite{KR07}.)
Since $A$ is $1$-Iwanaga-Gorenstein, this stably $3$-Calabi-Yau condition is equivalent to having functorial isomorphisms
$$
\Ext_A^1(M,N) \cong D\!\Ext_A^2(N,M) 
$$
for all $M \in \GP(A)$ and $N \in \md(A)$.

Namely, let $M \in \GP(A)$ and $N \in \md(A)$, and assume that $A$ is stably $3$-Calabi-Yau.
We have a short exact sequence
$$
0 \to \Omega_A(N) \to P(N) \to N \to 0
$$
with $P(N)$ projective.
Note that $\Omega_A(N) \in \GP(A)$ and $\Ext_A^1(M,P(N)) = 0$.
Applying $\Hom_A(M,-)$ and $\Hom_A(-,M)$ to this sequence and using the stably $3$-Calabi-Yau condition,
we get functorial isomorphisms
$$
\Ext_A^1(M,N) \cong \Ext_A^2(M,\Omega_A(N))
\cong D\!\Ext_A^1(\Omega_A(N),M)
\cong D\!\Ext_A^2(N,M).
$$

\begin{Thm}[Amiot, Keller--Reiten]\label{thm:Jacobi3CY}
Assume that $A = \Jacobi(Q,\tilde{W})$ is a Jacobian algebra defined by a hyperpotential $\tilde{W}$.
Then $A$ is $1$-Iwanaga-Gorenstein and stably
$3$-Calabi-Yau.
\end{Thm}

\begin{proof}
Let $A = \Jacobi_K(Q,\tilde{W})$ be a Jacobian algebra defined
by a hyperpotential $\tilde{W}$ for $Q$ in the sense of Ladkani \cite{L14}.
Then Amiot \cite[Corollary~3.6]{A09}  (for potentials)
and Ladkani \cite[Section~1.1]{L14} (for hyperpotentials) 
proved that $A$ is a $2$-Calabi-Yau tilted algebra.
Now \cite[Proposition~2.1(d) and Theorem~3.3]{KR07} implies that $A$ is $1$-Iwanaga-Gorenstein and stably
$3$-Calabi-Yau.
\end{proof}

%%%%%%%%%%%%%%%%%%%%%%%%%%%%%%%%%%%%%%%
%%%%%%%%%%%%%%%%%%%%%%%%%%%%%%%%%%%%%%%

\section{Standard sequences, $g$-vectors and $E$-symmetry}
\label{sec:2}

%%%%%%%%%%%%%%%%%%%%%%%%%%%%%%%%%%%%%%%
%%%%%%%%%%%%%%%%%%%%%%%%%%%%%%%%%%%%%%%

%%%%%%%%%%%%%%%%%%%%%%%%%%%%%%%%%%
\subsection{Dimension vectors}
%%%%%%%%%%%%%%%%%%%%%%%%%%%%%%%%%%
For $M \in \md(A)$, the Jordan-H\"older multiplicity of $S(i)$ 
for $M$ is denoted by $[M:S(i)]$.
The \emph{dimension vector} of $M$ is
$$
\dimv(M) := ([M:S(i)])_{i \in Q_0}.
$$

Let $q_{ij}$ be the number of arrows $i \to j$ in $Q$.

The following lemma is obvious.

\begin{Lem}\label{lem:dimv1}
We have
$$
\sum_{i,j \in Q_0} (q_{ji} - q_{ij}) [M:S(i)] [M:S(j)] = 0.
$$
\end{Lem}

Let $M \in \md(A)$.
For $i \in Q_0$, define $M(i) := e_iM$.
We identify 
$$
\Hom_A(P(i),M) \equiv M(i)
$$ via the map
$f \mapsto f(e_i)$.

For each $c \in e_iAe_j$ there is an $A$-module homomorphism
\begin{align*}
\cdot c \df P(i) &\to P(j)
\\
ae_i &\mapsto ae_ic.
\end{align*}
Applying $\Hom_A(-,M)$ to $\cdot c$ yields a $K$-linear map
\begin{align*}
M(c) := \Hom_A(\cdot c,M)\df M(j) &\to M(i)
\\
m &\mapsto cm.
\end{align*}
For $c \in e_iKQe_j$, we define $M(c) := M(\overline{c})$.

%%%%%%%%%%%%%%%%%%%%%%%%%%%%%%%%%%%%%%%%%%
\subsection{Standard complexes associated to modules}
%%%%%%%%%%%%%%%%%%%%%%%%%%%%%%%%%%%%%%%%%%
We recall the construction of standard complexes from \cite{GHLS26}.

An element $a \in KQ$ (resp. $a \in A$) is \emph{uniform}, if $a \in e_jKQe_i$ (resp. $a \in e_jAe_i$) for
some $i,j \in Q_0$.
In this case, we define $s(a) := i$ and $t(a) := j$.
We fix now a finite set
$
\cR
$ of uniform elements of $KQ$ such that $\cR$ generates the ideal $I$.

Recall, that in Section~\ref{subsec:1.2}, we defined $K$-linear maps
$$
l_a,r_a\df KQ \to KQ
\text{\quad and \quad}
\overline{l_a},\overline{r_a}\df KQ \to A.
$$

Note that $\Ima(l_a) \subseteq  e_{s(a)}KQ$ and 
$\Ima(r_a) \subseteq KQe_{t(a)}$.

For $\rho \in \cR$, we have
\[
\rho = \sum_{a \in Q_1 : \; t(a)=t(\rho)} a\,l_a(\rho)
= \sum_{b \in Q_1 : \; s(b)=s(\rho)} r_b(\rho)\,b.
\]
For a module $M$ and $k \in Q_0$, define
\[
\begin{aligned}
M_{\rm in}(k)
&:= \bigoplus_{a \in Q_1 : \; t(a)=k} M(s(a)),
&
M_{\rm out}(k)
&:= \bigoplus_{b \in Q_1 : \; s(b)=k} M(t(b)),\\
M'_{\rm rel}(k)
&:= \bigoplus_{\rho \in \cR : \; t(\rho)=k} M(s(\rho)),
&
M''_{\rm rel}(k)
&:= \bigoplus_{\rho \in \cR : \; s(\rho)=k} M(t(\rho)).
 \end{aligned}
\]
We get the \emph{standard complexes} 
\[
 M_{\rm rel}'(k)\xrightarrow{M(\gamma_{k,\cR}')}
 M_{\rm in}(k)\xrightarrow{M(\alpha_k)}M(k)\longrightarrow0,
\]
\[
 0\longrightarrow M(k)\xrightarrow{M(\beta_k)}
 M_{\rm out}(k)\xrightarrow{M(\gamma_{k,\cR}'')}M_{\rm rel}''(k),
\]
where
\begin{align*}
 M(\alpha_k)((v_a)_a)&:=\sum_{a \in Q_1 : \; t(a)=k}M(a)v_a,
 &M(\beta_k)(v)&:=(M(b)v)_{b},
\\
M(\gamma_{k,\cR}') &:=   (M(\gamma_{k,\cR}')_{a,\rho})_{a,\rho},
&
M(\gamma_{k,\cR}'') &:=   (M(\gamma_{k,\cR}'')_{\rho,b})_{\rho,b},
\end{align*}
where 
\begin{align*}
M(\gamma_{k,\cR}')_{a,\rho}\df M(s(\rho)) &\to M(s(a)),
&
M(\gamma_{k,\cR}'')_{\rho,b}\df M(t(b)) &\to M(t(\rho))
\\
w_\rho &\mapsto M(\overline{l_a}(\rho))w_\rho
&
v_b  &\mapsto M(\overline{r_b}(\rho))v_b.
\end{align*}

The following result from \cite[Section~5]{GHLS26} says that the standard complexes are indeed 
complexes of $K$-vector spaces.

\begin{Lem}\label{lem:complex1}
We have
$$
M(\alpha_k) \circ M(\gamma_{k,\cR}') = 0 
\text{\quad and \quad}
M(\gamma_{k,\cR}'') \circ M(\beta_k) = 0.
$$
\end{Lem}

The homology of the standard complexes yields a minimal
projective and a minimal injective presentation of $M$, see \cite[Theorem~1.1]{GHLS26} for details.

%%%%%%%%%%%%%%%%%%%%%%%%%%%%%%%%%%%%%%%
\subsection{Representation theoretic $g$-vectors}
%%%%%%%%%%%%%%%%%%%%%%%%%%%%%%%%%%%%%%%
For $M \in \md(A)$, let
\[
\bigoplus_{i \in Q_0} P(i)^{a_i} \to 
\bigoplus_{i \in Q_0} P(i)^{b_i} \to M \to 0
\]
be a minimal projective presentation.
Then the \emph{representation theoretic $g$-vector} $g^{\rm proj}(M)$ is defined as
$$
g^{\rm proj}(M) := (a_i-b_i)_{i \in Q_0}.
$$

Dually, for $M$ in $\md(A)$,
let
\[
0 \to M \to \bigoplus_{i \in Q_0} I(i)^{c_i} \to \bigoplus_{i \in Q_0} 
I(i)^{d_i} 
\]
be a minimal injective presentation. 
Then let 
$$
g^{\rm inj}(M) := (d_i-c_i)_{i \in Q_0}.
$$

For $M,N \in \md(A)$ and $i \ge 1$, let
$$
\ext_A^i(M,N) := \dim \Ext_A^i(M,N).
$$

\begin{Lem}\label{lem:gvector1}
For $M \in \md(A)$ and $i \in Q_0$, we have
\begin{align*}
g_i^{\rm proj}(M) &= \ext_A^1(M,S(i)) - \hom_A(M,S(i)),
\\
g_i^{\rm inj}(M) &= \ext_A^1(S(i),M) - \hom_A(S(i),M).
\end{align*}
\end{Lem}

\begin{proof}
We have 
\begin{align*}
a_i &= \ext_A^1(M,S(i)),  & b_i &= \hom_A(M,S(i)),
\\ 
d_i &= \ext_A^1(S(i),M), 
& c_i &= \hom_A(S(i),M).
\end{align*}
\end{proof}

%%%%%%%%%%%%%%%%%%%%%%%%%%%%%%%%%%%%%%%
\subsection{Rank defect and $E$-defect of modules}
%%%%%%%%%%%%%%%%%%%%%%%%%%%%%%%%%%%%%%%
For $i \in Q_0$,
let 
$$
\delta_i(M) := \rk(M(\gamma_{i,\cR}'')) - \rk(M(\gamma_{i,\cR}')) 
$$
be the \emph{rank defect} of $M$ at $i$.

\begin{Lem}\label{lem:gvector2}
For $M \in \md(A)$ and $i \in Q_0$, the following hold:
\begin{itemize}\itemsep2mm

\item[(i)]
$g_i^{\rm proj}(M) = \sum_{j \in Q_0} q_{ji} [M:S(j)]
- \rk(M(\gamma_{i,\cR}')) - [M:S(i)]$;

\item[(ii)]
$g_i^{\rm inj}(M) = \sum_{j \in Q_0} q_{ij} [M:S(j)] 
- \rk(M(\gamma_{i,\cR}'')) - [M:S(i)]$;

\item[(iii)]
$g_i^{\rm proj}(M) - g_i^{\rm inj}(M) = \sum_{j \in Q_0} 
(q_{ji} - q_{ij}) [M:S(j)] + \delta_i(M)$.

\end{itemize}
\end{Lem}

\begin{proof}
(i):
We have
\begin{align*}
g_i^{\rm proj}(M) &= \dim \Ker(M(\alpha_i)) - \rk(M(\gamma_{i,\cR}'))
- \dim \Coker(M(\alpha_i))
\\
&= (\dim M_{\rm in}(i) - \rk(M(\alpha_i))) - \rk(M(\gamma_{i,\cR}'))
- ([M:S(i)] - \rk(M(\alpha_i)))
\\
&= \dim M_{\rm in}(i) -  \rk(M(\gamma_{i,\cR}')) - [M:S(i)]
\\
&= \sum_{j \in Q_0} q_{ji} [M:S(j)]
- \rk(M(\gamma_{i,\cR}')) - [M:S(i)].
\end{align*}
Note that
the first equality is non-trivial and is part of  
\cite[Theorem~1.1]{GHLS26}.

(ii): This follows by duality.

(iii): This is a direct consequence of (i) and (ii).
\end{proof}

\begin{Lem}\label{lem:gvector3}
For $M,N \in \md(A)$ and $i,k \in Q_0$, 
the following hold:
\begin{itemize}\itemsep2mm

\item[(i)]
$\delta_i(M \oplus N) = \delta_i(M) + \delta_i(N)$;

\item[(ii)]
$\delta_i(M) = \delta_i(M \oplus S(k))$.

\end{itemize}
\end{Lem}

\begin{proof}
(i): 
This follows from the construction of the maps $M(\gamma_{i,\cR}')$
and $M(\gamma_{i,\cR}'')$ for $M \in \md(A)$.

(ii):
Since $I \subseteq \m^2$, $l_a(\m^2) \subseteq \m$,
and $r_a(\m^2) \subseteq \m$, we have 
$$
\rk(S(k)(\gamma_{i,\cR}')) = \rk(S(k)(\gamma_{i,\cR}'')) = 0.
$$
Now (ii) follows from (i).
\end{proof}

%%%%%%%%%%%%%%%%%%%%%%%%%%%%%%%%%%%%%%%
\subsection{$E$-invariants of finite-dimensional modules}
%%%%%%%%%%%%%%%%%%%%%%%%%%%%%%%%%%%%%%%
Recall that
$$
E(M) = \hom_A(M,\tau_A(M))
\text{\quad and \quad}
E^-(M) = \hom_A(\tau_A^-(M),M)
$$
are the \emph{$E$-invariants} of $M$.

For $M \in \md(A)$, let
$$
\partial(M) := E(M) - E^-(M)
$$
be the \emph{$E$-defect} of $M$.

For $a = (a_i)_i, b = (b_i)_i \in \Z^{Q_0}$, define
$$
\bil{a,b} := \sum_{i \in Q_0} a_ib_i.
$$
The following lemma is well known.

\begin{Lem}\label{lem:Einv1}
We have
\begin{align*}
E(M) &= \hom_A(M,M) + \bil{g^{\rm proj}(M),\dimv(M)}, 
\\
E^-(M) &= \hom_A(M,M) + \bil{g^{\rm inj}(M),\dimv(M)}.
\end{align*}
\end{Lem}

\begin{Cor}\label{cor:Einv1}
For $M \in \md(A)$, we have
$$
\partial(M) =  \sum_{i \in Q_0} [M:S(i)] 
(g_i^{\rm proj}(M) - g_i^{\rm inj}(M)).
$$
\end{Cor}

Following \cite{BS26},
the algebra $A$ is \emph{$E$-symmetric} if
$$
\partial(M) = 0
$$
for all $M \in \md(A)$.

\begin{Lem}\label{lem:Einv2}
For $M \in \md(A)$, we have
$$
\partial(M) = \sum_{i \in Q_0} [M:S(i)] \delta_i(M).
$$
\end{Lem}

\begin{proof}
This follows by combining Lemma~\ref{lem:gvector2}(iii) and Corollary~\ref{cor:Einv1}
 with Lemma~\ref{lem:dimv1}.
\end{proof}

\begin{Lem}\label{lem:Einv3}
For $M \in \md(A)$ and $k \in Q_0$, we have
$$
\partial(M \oplus S(k)) = \partial(M) + \delta_k(M).
$$
\end{Lem}

\begin{proof}
We have 
\begin{align*}
\partial(M \oplus S(k)) &= \sum_{i \in Q_0} [M \oplus S(k) :S(i)] \delta_i(M \oplus S(k))
\\
&= \sum_{i \in Q_0} [M:S(i)] \delta_i(M) + \delta_k(M)
\\
&= \partial(M) + \delta_k(M).
\end{align*}

The first and the third equality follows directly from Lemma~\ref{lem:Einv2}.
The second equality follows from 
Lemma~\ref{lem:gvector3}(ii).
\end{proof}

%%%%%%%%%%%%%%%%%%%%%%%%%%%%%%%%%%%%%%%%%%
\subsection{Known results on $E$-symmetric algebras}
%%%%%%%%%%%%%%%%%%%%%%%%%%%%%%%%%%%%%%%%%%

We recall some results from \cite{BS25} and \cite{BS26}.
The proofs work for arbitrary fields $K$.
(The assumption that $K$ is algebraically closed was only used for the geometric results in \cite{BS25} and \cite{BS26}.)

\begin{Thm}[{\cite[Theorem~1.7]{BS25}}]\label{thm:BS1}
Let $A$ be $E$-symmetric.
Then $A/AeA$ is $1$-Iwanaga-Gorenstein for all $e=e^2$ in $A$. In particular, we have
$\gldim(A) \in \{ 0, 1, \infty \}$.
\end{Thm}

\begin{Thm}[{\cite[Corollary~1.8]{BS25}}]\label{thm:BS2}
Assume that $\gldim(A) < \infty$.
Then the following are equivalent:
\begin{itemize}\itemsep2mm

\item[(i)]
$A$ is $E$-symmetric;

\item[(ii)]
$A$ is hereditary.

\end{itemize}
\end{Thm}

\begin{Prop}[{\cite[Proposition~3.2]{BS26}}]\label{prop:BS3}
The following are equivalent:
\begin{itemize}\itemsep2mm

\item[(i)]
$E(M,N) = E^-(M,N)$ for all $M,N \in \md(A)$;

\item[(ii)]
$A$ is hereditary.

\end{itemize}
\end{Prop}
%%%%%%%%%%%%%%%%%%%%%%%%%%%%%%%%%%%%%%%%%%
\subsection{Characterization of $E$-symmetry via ranks}
%%%%%%%%%%%%%%%%%%%%%%%%%%%%%%%%%%%%%%%%%%

\begin{Prop}\label{prop:EsymmRank}
The following are equivalent:
\begin{itemize}\itemsep2mm

\item[(i)]
$A$ is $E$-symmetric;

\item[(ii)]
For all $M \in \md(A)$ and $i \in Q_0$, we have
$$
\rk(M(\gamma_{i,\cR}')) = \rk(M(\gamma_{i,\cR}''));
$$

\item[(iii)]
For all $M \in \md(A)$ and $i \in Q_0$, we have
$$
g_i^{\rm proj}(M) - g_i^{\rm inj}(M)
= \sum_{j \in Q_0} (q_{ji} - q_{ij})[M:S(j)].
$$

\end{itemize}
\end{Prop}

\begin{proof}
(i) $\implies$ (ii):
Assuming (i), we get $\partial(M \oplus S(k)) = \partial(M) = 0$.
Lemma~\ref{lem:Einv3} implies now $\delta_k(M) = 0$ for all
$k \in Q_0$.
This proves (ii).

(ii) $\implies$ (iii):
Condition (ii) says that $\delta_i(M) = 0$ for all $i \in Q_0$.
Now Lemma~\ref{lem:gvector2}(iii) implies (iii).

(iii) $\implies$ (i):
Assuming (iii), we get
\begin{align*}
\partial(M) &= \sum_{i \in Q_0} [M:S(i)] 
(g_i^{\rm proj}(M) - g_i^{\rm inj}(M))
\\
&= \sum_{i,j \in Q_0} [M:S(i)](q_{ji} - q_{ij})[M:S(j)]
\\
&= 0.
\end{align*}
The first equality is Corollary~\ref{cor:Einv1}, and
for the third equality, we used Lemma~\ref{lem:dimv1}.
Thus $A$ is $E$-symmetric.
\end{proof}

%%%%%%%%%%%%%%%%%%%%%%%%%%%%%%%%%%%%%%
\subsection{Hom criterion for isomorphic modules}
%%%%%%%%%%%%%%%%%%%%%%%%%%%%%%%%%%%%%%
The following theorem is due to Auslander. A proof 
(of a more general version) can be found in \cite{B89}.

\begin{Thm}[Auslander]\label{thm:auslander}
Let $M,N \in \md(A)$. Then the following hold:
\begin{itemize}\itemsep2mm

\item[(i)]
If $\hom_A(M,X) = \hom_A(N,X)$ for all $X \in \md(A)$,
then $M \cong N$;

\item[(ii)]
If $\hom_A(X,M) = \hom_A(X,N)$ for all $X \in \md(A)$,
then $M \cong N$.

\end{itemize}
\end{Thm}

\begin{Lem}\label{lem:cokrank}
For $f,g \in \Hom_A(P_1,P_0)$, the following are equivalent:
\begin{itemize}\itemsep2mm

\item[(i)]
There exist $h_1 \in \Aut_A(P_1)$ and $h_0 \in \Aut_A(P_0)$ with
$$
g = h_0^{-1} \circ f \circ h_1;
$$

\item[(ii)]
$\Coker(f) \cong \Coker(g)$;

\item[(iii)]
For each $M \in \md(A)$, we have
$$
\rk(\Hom_A(f,M)) = \rk(\Hom_A(g,M)).
$$ 

\end{itemize}
\end{Lem}

\begin{proof}
Since $\Hom_A(-,M)$ is left exact, we get
\begin{equation}\label{eq:6}
\rk(\Hom_A(h,M)) = \hom_A(P_0,M) - \hom_A(\Coker(h),M)
\end{equation}
for $h \in \{ f,g \}$.

(i) $\iff$ (ii): This is straightforward and well known.

(ii) $\implies$ (iii):
Condition (ii) together with \eqref{eq:6} imply (iii).

(iii) $\implies$ (i):
Condition (iii) combined with \eqref{eq:6} imply that
$$\hom_A(\Coker(f),M) = \hom_A(\Coker(g),M)$$ for
all $M \in \md(A)$.
Now Theorem~\ref{thm:auslander} implies $\Coker(f) \cong \Coker(g)$.
\end{proof}

%%%%%%%%%%%%%%%%%%%%%%%%%%%%%%%%%%%%%%%
%%%%%%%%%%%%%%%%%%%%%%%%%%%%%%%%%%%%%%%

\section{Characterization of $E$-symmetric algebras via simple modules}\label{sec:3}

%%%%%%%%%%%%%%%%%%%%%%%%%%%%%%%%%%%%%%%
%%%%%%%%%%%%%%%%%%%%%%%%%%%%%%%%%%%%%%%

For $M \in \md(A)$, we denote its projective cover (resp. injective envelope)
by
$$
P(M) \xrightarrow{p_M} M \to 0
\qquad
(\text{resp.} \quad 0 \to M \xrightarrow{i_M} I(M))
$$

For $i \in Q_0$, we get short exact sequences
\begin{equation}\label{eq:1a}
0 \to \Omega_A(S(i)) \to P(i) \xrightarrow{p_{S(i)}} S(i) \to 0
\end{equation}
and 
\begin{equation}\label{eq:1b}
0 \to S(i) \xrightarrow{i_{S(i)}} I(i) \to \Omega_A^{-1}(S(i))
\to 0.
\end{equation}
We have
$$
\Omega_A(S(i)) = \rad(P(i))
\text{\quad and \quad}
\Omega_A^{-1}(S(i)) = I(i)/S(i).
$$
Let 
$$
I(i)/S(i) = X(i) \oplus \bigoplus_{j \in Q_0} I(j)^{c_{ij}},
$$
where $X(i)$ has no non-zero injective direct summand.

Recall that $q_{ij}$ is the number of arrows $i \to j$ in $Q_1$.
We have
$$
q_{ij} = \ext_A^1(S(i),S(j)).
$$

Furthermore, let
$$
r_{ij} := \ext_A^2(S(i),S(j)).
$$

\begin{Thm}\label{thm:mainpart1}
The following are equivalent:
\begin{itemize}\itemsep2mm

\item[(i)]
$A$ is $E$-symmetric;

\item[(ii)]
In $\overline{\md}(A)$ we have
$\tau_A(\Omega_A(S(i))) \cong I(i)/S(i)$ for each $i \in Q_0$;

\item[(iii)]
We have $\tau_A(\Omega_A(S(i))) \cong X(i)$ for each $i \in Q_0$;

\item[(iv)]
$A$ is a $1$-Iwanaga-Gorenstein algebra, and for each simple $A$-module $S$ and $M \in \GP(A)$, we have an isomorphism
$$
\Ext_A^1(M,S) \cong D\!\Ext_A^2(S,M)
$$
which is functorial in $M$.

\end{itemize}
In this case, we have $r_{ij} + c_{ij} = q_{ji}$.
\end{Thm}

\begin{proof}
(iii) $\implies$ (ii): Obvious.

(ii) $\implies$ (iii):
The module $\tau_A(\rad(P_i))$ has no non-zero injective direct summand.
The claim follows.

(i) $\iff$ (iii):

Let $M \in \md(A)$.
We have $[M:S(i)] = \hom_A(P(i),M) = \hom_A(M,I(i))$.

Applying $\Hom_A(M,-)$ to the short exact sequence \eqref{eq:1b} and $\Hom_A(-,M)$ to 
\eqref{eq:1a}  and using Lemma~\ref{lem:gvector1} yields
\begin{align*}
g_i^{\rm proj}(M) &= \hom_A(M,I(i)/S(i)) - [M:S(i)],
\\
g_i^{\rm inj}(M) &= \hom_A(\rad(P(i)),M) - [M:S(i)].
\end{align*}
The condition Proposition~\ref{prop:EsymmRank}(iii) is then equivalent to
\begin{equation}\label{eq:2}
\hom_A(M,I(i)/S(i)) - \hom_A(\rad(P(i)),M) =
\sum_{j \in Q_0} (q_{ji}-q_{ij}) [M:S(j)].
\end{equation}

Applying $\Hom_A(-,S(j))$ to \eqref{eq:1a} yields that 
the minimal projective presentation of $\rad(P(i))$ has the form
\begin{equation}\label{eq:rqprojpresent}
\bigoplus_{j \in Q_0} P(j)^{r_{ij}} \to 
\bigoplus_{j \in Q_0} P(j)^{q_{ij}} \to \rad(P(i)) \to 0.
\end{equation}
Now \cite[Lemma~2.7]{BS25} implies that
\begin{equation}\label{eq:3}
\hom_A(\rad(P(i)),M) - \hom_A(M,\tau_A(\rad(P(i)))) = 
\sum_{j \in Q_0} (q_{ij} - r_{ij}) [M:S(j)].
\end{equation}
Now we insert \eqref{eq:3} into \eqref{eq:2}.
Simplifying the resulting equation gives that 
Proposition~\ref{prop:EsymmRank}(iii) is equivalent to
\begin{equation}\label{eq:4}
\hom_A(M,I(i)/S(i)) - \hom_A(M,\tau_A(\rad(P(i)))) = 
\sum_j(q_{ji}-r_{ij}) [M:S(j)].
\end{equation}
Since $\hom_A(M,I(j)) = [M:S(j)]$, \eqref{eq:4} is equivalent 
to
$$
\hom_A\left(M,\tau_A(\rad(P(i))) \oplus \bigoplus_{j \in Q_0} 
I(j)^{q_{ji}}\right) =
\hom_A\left(M, I(i)/S(i) \oplus\bigoplus_{j \in Q_0} I(j)^{r_{ij}}\right).
$$
Since this equality holds now for all $M \in \md(A)$, we
can apply 
Theorem~\ref{thm:auslander} and get that \eqref{eq:4} is equivalent to
\begin{equation}\label{eq:5}
 \tau_A(\rad(P(i))) \oplus \bigoplus_{j \in Q_0} I(j)^{q_{ji}}
       \cong I(i)/S(i) \oplus\bigoplus_{j \in Q_0} I(j)^{r_{ij}}.
\end{equation}
Equation~\eqref{eq:5} together with the Krull-Remak-Schmidt Theorem implies that $q_{ji} = c_{ij}+r_{ij}$ and $\tau_A(\rad(P(i))) \cong X(i)$.
Vice versa, assume that $\tau_A(\rad(P(i))) \cong X(i)$.
Applying the Nakayama functor $\nu_A$ to the 
minimal projective presentation \eqref{eq:rqprojpresent} yields a minimal
injective presentation of $\tau_A(\rad(P(i)))$.
In particular, we get
$$
\soc(\tau_A(\rad(P(i)))) \cong \bigoplus_{j \in Q_0} S(j)^{r_{ij}}.
$$
Furthermore, we have
$$
\soc(I(i)/S(i)) \cong \bigoplus_{j \in Q_0} S(j)^{q_{ji}}
\text{\quad and \quad}
\soc(X(i)) \cong \bigoplus_{j \in Q_0} S(j)^{q_{ji}-c_{ij}}.
$$
This implies
$q_{ji} = c_{ij}+r_{ij}$, and therefore \eqref{eq:5} holds.

(ii) $\implies$ (iv):
Assume (ii) holds.

Claim:
For each $M \in \GP(A)$, we have 
functorial isomorphisms
\begin{align*}
D\!\Ext_A^2(S,M) 
&\cong D\!\Ext_A^1(\Omega_A(S),M)
\cong \overline{\Hom}_A(M,\tau_A(\Omega_A(S)))
\\
&\cong \overline{\Hom}_A(M,\Omega_A^{-1}(S)) 
\cong \Ext_A^1(M,S).
\end{align*}
In particular, (iv) holds.

Proof:
The $1$st isomorphism is obvious, the $2$nd holds by the Auslander-Reiten formulas, and the $3$rd follows from assumption (ii).
We have to prove the $4$th isomorphism:
Applying $\Hom_A(M,-)$ to an exact sequence
$$
0 \to S \xrightarrow{i_S} I(S) \xrightarrow{q} \Omega_A^{-1}(S) \to 0,
$$
we get a functorial isomorphism
$$
\Ext_A^1(M,S) \cong \Hom_A(M,\Omega_A^{-1}(S))/\Ima(\Hom_A(M,q)).
$$
(Note that $\Ext_A^1(M,I(S)) = 0$, since $I(S)$ is injective.)
By construction, all homomorphisms in $\Ima(\Hom_A(M,q))$ factor through the
injective module $I(S)$.
Let now $f \in \Hom_A(M,\Omega_A^{-1}(S))$ such that $f = u_2u_1$, where
$u_1 \in \Hom_A(M,I)$ and $u_2 \in \Hom_A(I,\Omega_A^{-1}(S))$ with $I$ injective.
Since $A$ is $1$-Iwanaga-Gorenstein, we have $\pdim(I) \le 1$.
Thus there is a projective resolution
$$
0 \to P_1 \to P_0 \xrightarrow{p_I} I \to 0.
$$
Since $M \in \GP(A)$, we have $\Ext_A^1(M,P_1) = 0$.
Applying now $\Hom_A(M,-)$ to this resolution of $I$, we get that
there exists some $u\df M \to P_0$ with $u_1 = p_Iu$.
The projectivity of $P_0$ gives a homomorphism $w\df P_0 \to I(S)$ such that
$u_2p_I = qw$.
The following diagram illustrates the situation:
$$
\xymatrix{
&&& M \ar[d]^{u_1}\ar@{-->}[dl]_{u}
\\
0 \ar[r] & P_1 \ar[r] & P_0 \ar[r]^{p_I}\ar@{-->}[d]^w & I \ar[r]\ar[d]^{u_2} & 0
\\
0 \ar[r] & S \ar[r]^{i_S} & I(S) \ar[r]^q& \Omega_A^{-1}(S) \ar[r] & 0
}
$$
We get
$$
f = u_2u_1 = u_2p_Iu = qwu.
$$
Thus $f$ factors through $q$, and we conclude that
$$
\Ext_A^1(M,S) \cong \overline{\Hom}_A(M,\Omega_A^{-1}(S)).
$$
This finishes the proof of the claim.
$\square$

(iv) $\implies$ (ii):
Assume (iv) holds.

Claim:
For each $X \in \md(A)$, we have isomorphisms
\begin{align*}
\underline{\Hom}_A(\Omega_A(S),X) &\cong \Ext_A^1(\Omega_A(S),\Omega_A(X))
\cong \Ext_A^2(S,\Omega_A(X))
\\
&\cong D\!\Ext_A^1(\Omega_A(X),S)
\cong D\!\Ext_A^2(X,S)
\\
&\cong D\!\Ext_A^1(X,\Omega_A^{-1}(S))
\cong \underline{\Hom}_A(\tau_A^{-1}(\Omega_A^{-1}(S)),X)
\end{align*}
which are functorial in $X$.

Proof:
The 
$3$rd isomorphism comes from assumption (iv). 
The $2$rd, $4$th and $5$th are obvious, and the $6$th holds by the Auslander-Reiten formulas. 
We now prove the $1$st isomorphism:
There is a short exact sequence
$$
0 \to \Omega_A(X) \to P(X) \xrightarrow{p_X} X \to 0
$$
with $P(X)$ projective.
We have $\Ext_A^1(\Omega_A(S),P(X)) = 0$, since $\Omega_A(S) \in \GP(A)$.
Thus, applying $\Hom_A(\Omega_A(S),-)$ to the exact sequence above, we get
an isomorphism
\begin{align*}
\Ext_A^1(\Omega_A(S),\Omega_A(X)) &\cong \Hom_A(\Omega_A(S),X)/\Ima(\Hom_A(\Omega_A(S),p_X))
\\
&= \underline{\Hom}_A(\Omega_A(S),X),
\end{align*}
which is functorial in $X$.
$\square$

The resulting isomorphism 
$$
\underline{\Hom}_A(\Omega_A(S),-) \cong \underline{\Hom}_A(\tau_A^{-1}(\Omega_A^{-1}(S)),-)
$$
of functors
implies by the Yoneda Lemma that 
$$
\Omega_A(S) \cong \tau_A^{-1}(\Omega_A^{-1}(S))
$$ in $\underline{\md}(A)$.
(This is the place where we use the functoriality in condition (iv), otherwise we could not apply the Yoneda Lemma.)
This is equivalent to 
$$
\tau_A(\Omega_A(S)) \cong \Omega_A^{-1}(S)
$$ in $\overline{\md}(A)$.
Thus (ii) holds.
\end{proof}

%%%%%%%%%%%%%%%%%%%%%%%%%%%%%%%%%%%%%%%
%%%%%%%%%%%%%%%%%%%%%%%%%%%%%%%%%%%%%%%

\section{$E$-symmetric algebras and generalized 
Jacobian algebras}\label{sec:4}

%%%%%%%%%%%%%%%%%%%%%%%%%%%%%%%%%%%%%%%
%%%%%%%%%%%%%%%%%%%%%%%%%%%%%%%%%%%%%%%

%%%%%%%%%%%%%%%%%%%%%%%%%%%%%%%%%%%%%%%%%
\subsection{The relation space ideal}
\label{subsec:9}
%%%%%%%%%%%%%%%%%%%%%%%%%%%%%%%%%%%%%%%%%
For $i \ge 0$, let $Q_i$ be the set of paths of length $i$ in $Q$,
and let $KQ_i$ be the $K$-span of these paths.
 Let
$$
\m := \bigoplus_{i \ge 1} KQ_i 
$$
be the \emph{arrow ideal} of $KQ$.
We have
$$
\m^{N_I} \subseteq I \subseteq \m^2
$$
for some $N_I \ge 2$.
Then $J(A) := \m/I$ is the Jacobson radical of $A$.

As before, 
for $i,j \in Q_0$, let
$$
r_{ij} := \ext_A^2(S(i),S(j)).
$$

We fix now some
$N_I \ge 2$ such that $\m^{N_I} \subseteq I$.

We call 
$$
R := \frac{I}{\m I+I\m}
$$
the \emph{relation space} for $A$.

\begin{Lem}\label{lem:truncation1}
We have
$$
\dim e_j\frac{I}{\m I+I\m}e_i = r_{ij}.
$$
\end{Lem}

\begin{proof}
This can be found in \cite[Section~1]{B83}.
(Bongartz assumes that the field $K$ is algebraically closed.
But this assumption is not needed for algebras of the form 
$A = KQ/I$.)
\end{proof}

By Lemma~\ref{lem:truncation1}, we can choose for each $i,j \in Q_0$ uniform elements
$$
\rho_{ij,l} \in e_jIe_i 
$$ with $1 \le l \le r_{ij}$ whose 
residue classes $\overline{\rho}_{ij,l}$ %(module $\m I + I \m$)
yield a $K$-basis 
$\{ e_j\overline{\rho}_{ij,l}e_i \mid 1 \le l \le r_{ij} \}$
of
$$
e_j \frac{I}{\m I + I \m} e_i.
$$

Let 
$$
W := \{ \rho_{ij,l} \mid i,j \in Q_0,\; 1 \le l \le r_{ij} \},
$$
and let $I(W)$ be the ideal of $KQ$ generated by $W$.
We call $I(W)$ the \emph{relation space ideal} of $A$.
One can assume without loss of generality that
$$
\rho_{ij,l} \in KQ_{\le N_I},
$$
since $\m^{N_I+1} \subseteq \m I + I \m$.
This does not affect the ideal $I(W)$.

\begin{Lem}\label{lem:truncation2}
We have
$$
I = I(W) + \m^{N_I+1}.
$$
\end{Lem}

\begin{proof} 
Since the set 
$$
\{ \overline{\rho}_{ij,l} \mid i,j \in Q_0,\; 1 \le l \le r_{ij} \}
$$
is a basis of the $K$-vector space
$$
\frac{I}{\m I + I \m},
$$
we have
\begin{equation}\label{eq:7}
I = I(W) + \m I + I \m.
\end{equation}
Now we can replace each $I$ appearing in the right hand side of \eqref{eq:7}
by the right hand side of \eqref{eq:7}
and obtain
\begin{align*}
I &= I(W) + \m(I(W) + \m I + I \m) + (I(W) + \m I + I \m)\m
\\
&= I(W) + \m^2 I + \m I \m + I \m^2. 
\end{align*}
Iterating this, we get by induction that
$$
I = I(W) +\sum_{d_1+d_2 = d} \m^{d_1} I \m^{d_2} 
$$
for each $d \ge 1$.

For $d = d_1+d_2 = N_I-1$ we use that $I \subseteq \m^2$ and get
$$
\m^{d_1} I \m^{d_2} \subseteq \m^{d_1+d_2+2} = \m^{N_I+1}.
$$
Thus $I \subseteq I(W) + \m^{N_I+1}$.
Conversely, we have
$I(W) + \m^{N_I+1} \subseteq I$, since $I(W) \subseteq I$ and
$\m^{N_I} \subseteq I$.
This finishes the proof.
\end{proof}

%%%%%%%%%%%%%%%%%%%%%%%%%%%%%%%%%%%%%%%%%%%%%%%%
\subsection{From $E$-symmetry to generalized hyperpotentials}
\label{subsec:h2}
%%%%%%%%%%%%%%%%%%%%%%%%%%%%%%%%%%%%%%%%%%%%%%%%
Assume now that $A$ is $E$-symmetric.
Then
Theorem~\ref{thm:mainpart1} implies that
$$
r_{ij} \le q_{ji}
$$
for all $i,j \in Q_0$.
Thus we can choose an
an injective map
$$
\eta_{ij}\df \{ \rho_{ij,l} \mid 1 \le l \le r_{ij} \}
\to \{ b \in Q_1 \mid s(b) = j,\; t(b) = i \}.
$$
For each $a = \eta_{ij}(\rho_{ij,l}) \in \Ima(\eta_{ij})$, let
$$
\rho_a := \rho_{ij,l} \in e_{s(a)} I e_{t(a)}.
$$
Here is the corresponding visualization:
$$
\xymatrix{
j \ar[r]^a 
& i \ar@{-->}@/^4ex/[l]^{\rho_a}
}
$$
For all other arrows in $a \in \{ b \in Q_1 \mid s(b) = j,\; t(b) = i \}
\setminus \Ima(\eta_{ij})$, define $\rho_a := 0$.

Thus, the relation space ideal $I(W)$ is generated by
$$
W := \{ \rho_a \mid a \in Q_1 \}.
$$
By Lemma~\ref{lem:truncation2},
$W$ satisfies condition (h1) in Definition~\ref{def:genhyper}.

We need to show that $W$ also satisfies (h2).

From here on, we work with the set 
$$
\cR := W^+ := \{ \rho_a \mid a \in Q_1 \} \cup Q_{N_I+1}
$$
of relations generating the ideal $I$.

We use the same notation as in Section~\ref{subsec:1.2}.

%We have
%$\scR_k, \scL_k \in 
%\rad_A(P_{k,{\rm in}},P_{k,{\rm out}})$.

Recall that for all $a,b \in Q_1$, we have 
$$
\rho_a \in e_{s(a)} I e_{t(a)}
\text{\quad and \quad}
\overline{r_b}(\rho_a),\, \overline{l_a}(\rho_b)
\in e_{s(a)}Ae_{t(b)}.
$$

Define
\begin{align*}
M_{N_I}'(k)&:= \bigoplus_{p \in Q_{N_I+1} :\; t(p) = k} M(s(p)),
&
M_{N_I}''(k) &:= \bigoplus_{p \in Q_{N_I+1} :\; s(p) = k} M(t(p)).
\end{align*}
Recall that

\begin{align*}
M_{\rm in}(k)
&:= \bigoplus_{a \in Q_1 : \; t(a)=k} M(s(a)),
&
M_{\rm out}(k)
&:= \bigoplus_{b \in Q_1 : \; s(b)=k} M(t(b)).
\end{align*}
We get
\begin{align*}
M_{\rm rel}'(k)
:= \bigoplus_{\rho \in W^+ : \; t(\rho)=k} M(s(\rho))
&= \bigoplus_{b \in Q_1 :\; s(b) = k} M(t(b)) 
\oplus M_{N_I}'(k)
\\
&= M_{\rm out}(k) \oplus M_{N_I}'(k),
\\
M_{\rm rel}''(k)
:= \bigoplus_{\rho \in W^+ : \; s(\rho)=k} M(t(\rho))
&= \bigoplus_{a \in Q_1 : \; t(a)=k} M(s(a)) \oplus M_{N_I}''(k)
\\
&= M_{\rm in}(k) \oplus M_{N_I}''(k).
\end{align*}

Recall that for each $i \in Q_0$, we identify
$$
\Hom_A(P(i),M) \equiv M(i)
$$ via $f \mapsto f(e_i)$.
This leads to the identifications
$$
\Hom_A(P_{k,{\rm out}},M) \equiv M_{\rm out}(k),
\text{\quad and \quad}
\Hom_A(P_{k,{\rm in}},M) \equiv M_{\rm in}(k).
$$

We can think of the maps 
\begin{align*}
M(\mu_k') := \Hom_A(\scL_k,M) = (\Hom_A(\scL_{k,b,a},M))_{a,b}
\df M_{\rm out}(k) &\to M_{\rm in}(k),
\\
M(\mu_k'') := \Hom_A(\scR_k,M)= (\Hom_A(\scR_{k,b,a},M))_{a,b}\df 
M_{\rm out}(k) &\to M_{\rm in}(k)
\end{align*}
as matrices whose entries with row index $a$ and column index 
$b$ are the maps $\Hom_A(\scL_{k,b,a},M)$ and 
$\Hom_A(\scR_{k,b,a},M)$.
It follows now directly from the definitions that
the following two
diagrams
$$
\xymatrix@+3ex{
\Hom_A(P(t(b)),M) \ar@3{-}[d] \ar[rr]^{\Hom_A(\scL_{k,b,a},M)} 
&& \Hom_A(P(s(a)),M) \ar@3{-}[d]
\\
M(t(b)) \ar[rr]^{M(\ov{l_a}(\rho_b))} && M(s(a))
}
$$
$$
\xymatrix@+3ex{
\Hom_A(P(t(b)),M) \ar@3{-}[d] \ar[rr]^{\Hom_A(\scR_{k,b,a},M)} 
&& \Hom_A(P(s(a)),M) \ar@3{-}[d]
\\
M(t(b)) \ar[rr]^{M(\ov{r_b}(\rho_a))} && M(s(a))
}
$$
are commutative.

\begin{Lem}\label{lem:h2}
For $k \in Q_0$ and $\cR = W^+$ as defined above, we have
\begin{align*}
M(\gamma_{k,\cR}') &= (M(\mu_k'),0)\df 
M_{\rm out}(k) \oplus M_{N_I}'(k) \to M_{\rm in}(k),
\\
M(\gamma_{k,\cR}'') &= \left(\begin{matrix} M(\mu_k'')\\0\end{matrix}\right)\df 
M_{\rm out}(k)  \to M_{\rm in}(k) \oplus M_{N_I}''(k).
\end{align*}
\end{Lem}

\begin{proof}
For each $p \in Q_{N_I+1}$ we have
$M(\overline{r_b}(p)) = M(\overline{l_a}(p)) = 0$ for all
$a,b \in Q_1$, since $r_b(p)$ and $l_a(p)$ are either $0$ or of length $N_I$.
But all paths of length $N_I$ annihilate all $A$-modules, since
$\m^{N_I} \subseteq I$.

Now the definitions of $M(\gamma_{k,\cR}')$ and 
$M(\gamma_{k,\cR}'')$ together with the preceeding explanations yield the result.
\end{proof}

\begin{Cor}\label{cor:h2}
For $k \in Q_0$ and $\cR = W^+$ as defined above, we have
\begin{align*}
\rk(M(\gamma_{k,\cR}')) &= \rk(\Hom_A(\scL_k,M)),
\\
\rk(M(\gamma_{k,\cR}'')) &= \rk(\Hom_A(\scR_k,M)).
\end{align*}
\end{Cor}

By Proposition~\ref{prop:EsymmRank}, the $E$-symmetry of $A$ implies now that
$\rk(\Hom_A(\scL_k,M)) = \rk(\Hom_A(\scR_k,M))$.
Now Lemma~\ref{lem:cokrank} yields automorphisms
\[
h_{k,{\rm in}} \in \Aut_A(P_{k,{\rm in}})
\text{\quad and \quad}
h_{k,{\rm out}} \in \Aut_A(P_{k,{\rm out}})
\]
such that
$$
\scL_k = h_{k,{\rm out}}^{-1} \circ \scR_k \circ h_{k,{\rm in}}.
$$
In other words, condition (h2) holds.
Thus $A$ is a generalized Jacobian algebra.

%%%%%%%%%%%%%%%%%%%%%%%%%%%%%%%%%%%%%%%%%%
\subsection{Characterization of $E$-symmetric algebras}
%%%%%%%%%%%%%%%%%%%%%%%%%%%%%%%%%%%%%%%%%%
We can now finish the proof of our main result:

\begin{Thm}\label{thm:mainpart2}
The following are equivalent:
\begin{itemize}\itemsep2mm

\item[(i)] 
$A$ is $E$-symmetric;

\item[(ii)] $A$ is generalized Jacobian in the sense of Definitions~\ref{def:genhyper} and \ref{def:genjacobian}.

\end{itemize}
\end{Thm}

\begin{proof}
(i) $\implies$ (ii):
We proved this implication in Section~\ref{subsec:h2}.

(ii) $\implies$ (i):
Under assumption (ii) we get the equation
$$
\scL_k = h_{k,{\rm out}}^{-1} \circ \scR_k \circ h_{k,{\rm in}}
$$
for some $h_{k,{\rm in}} \in \Aut_A(P_{k,{\rm in}})$ and $h_{k,{\rm out}} \in \Aut_A(P_{k,{\rm out}})$.
Now Lemma~\ref{lem:cokrank} implies that
$$
\rk(\Hom_A(\scL_k,M)) = \rk(\Hom_A(\scR_k,M)).
$$
Then
Proposition~\ref{prop:EsymmRank} together with
Corollary~\ref{cor:h2} implies that $A$ is $E$-symmetric.
\end{proof}

%%%%%%%%%%%%%%%%%%%%%%%%%%%%%%%%%%%%%%%%%%
%%%%%%%%%%%%%%%%%%%%%%%%%%%%%%%%%%%%%%%%%%

\section{Characterization of selfinjective $E$-symmetric algebras}\label{sec:5}

%%%%%%%%%%%%%%%%%%%%%%%%%%%%%%%%%%%%%%%%%%
%%%%%%%%%%%%%%%%%%%%%%%%%%%%%%%%%%%%%%%%%%

Assume that $A$ is selfinjective.
The syzygy functor $\Omega_A$ and the
Nakayama functor
$\nu_A := D\!\Hom_A(-,{_A}A)$
are 
then autoequivalences
$$
\Omega_A,\; \nu_A\df \underline{\md}(A) \to \underline{\md}(A)
$$
of the stable module category.
On $\underline{\md}(A)$, we have functorial isomorphisms
$$
\tau_A \cong \Omega_A^2 \circ \nu_A 
\text{\quad and \quad}
\Omega_A \circ \nu_A \cong \nu_A \circ \Omega_A,
$$
see for example \cite[Chapter~IV]{ARS97}. 

As a consequence of Theorem~\ref{thm:main1}, we obtain the
following characterization of $E$-symmetric selfinjective algebras.

\begin{Cor}\label{cor:main2b}
Let $A$ be a selfinjective algebra.
Then the following are equivalent:
\begin{itemize}\itemsep2mm

\item[(i)]
$A$ is $E$-symmetric;

\item[(ii)]
For each non-projective simple $A$-module $S$, we have
$$
\Omega_A^4(S) \cong \nu_A^{-1}(S).
$$
\end{itemize}
\end{Cor}

Note that Corollary~\ref{cor:main2b}(ii) is equivalent to
$$
\tau_A^2(S) \cong \nu_A(S)
$$
for each non-projective simple $A$-module $S$.

\begin{proof}
By 
Theorem~\ref{thm:main1}, $A$ is $E$-symmetric if and only if 
$$
\tau_A(\Omega_A(S)) \cong \Omega_A^{-1}(S)
$$ 
in $\overline{\md}(A)$ for each simple $A$-module $S$.
Since $A$ is selfinjective, we have
$\overline{\md}(A) = \underline{\md}(A)$ and
$\tau_A \cong \Omega_A^2 \circ \nu_A$ as autoequivalences of $\underline{\md}(A)$.
Thus the condition becomes equivalent to 
$$
\Omega^4_A(S) \cong \nu_A^{-1}(S)
$$
in $\underline{\md}(A)$ for each simple $A$-module $S$.
This is equivalent to having an isomorphism
$$
\Omega^4_A(S) \cong \nu_A^{-1}(S)
$$
in $\md(A)$ for each non-projective simple $A$-module $S$.
This finishes the proof.
\end{proof}

\begin{Rem}
For $E$-symmetric selfinjective algebras, the inequality
$r_{ij} \le q_{ji}$ resulting from Theorem~\ref{thm:mainpart1}
is an equality.
\end{Rem}

\begin{Cor}\label{cor:main4b}
Let $A$ be local.
Then the following are equivalent:
\begin{itemize}\itemsep2mm

\item[(i)]
$A$ is $E$-symmetric;

\item[(ii)]
$A$ is selfinjective, and if the simple $A$-module $S$ is
non-projective, then
$$
\Omega_A^4(S) \cong S.
$$

\end{itemize}
\end{Cor}

\begin{proof}
If $A$ is local and selfinjective and if the simple $A$-module $S$ is
non-projective, then 
$\nu_A^{-1}(S) \cong S$.
We will use this twice in the proof.

(i) $\implies$ (ii):
Let $A$ be $E$-symmetric.
Then \cite[Theorem~1.7]{BS25} implies that $A$ is 
$1$-Iwanaga-Gorenstein. But a local $1$-Iwanaga-Gorenstein algebra has to be selfinjective.
(Here we use that for a local algebra $A$, the $A$-modules $M$ with 
$\pdim(M) < \infty$ are projective.)
Now Corollary~\ref{cor:main2b} implies (ii).

(ii) $\implies$ (i): 
This follows directly from 
Corollary~\ref{cor:main2b}.
\end{proof}

%%%%%%%%%%%%%%%%%%%%%%%%%%%%%%%%%%%%%%%%%%
%%%%%%%%%%%%%%%%%%%%%%%%%%%%%%%%%%%%%%%%%%

\section{An $E$-symmetric algebra which is not a Jacobian 
algebra}\label{sec:6}

%%%%%%%%%%%%%%%%%%%%%%%%%%%%%%%%%%%%%%%%%%
%%%%%%%%%%%%%%%%%%%%%%%%%%%%%%%%%%%%%%%%%%

%%%%%%%%%%%%%%%%%%%%%%%%
\subsection{Enomoto's example is a generalized Jacobian algebra}
%%%%%%%%%%%%%%%%%%%%%%%%
Using the AI model ChatGPT-6 Astra, 
Enomoto \cite{E26} found a counterexample to the Periodicity Conjecture.
This example also serves as an example of a generalized Jacobian algebra which is not a Jacobian algebra defined by a hyperpotential (in the sense of Ladkani \cite{L14}).

Note that Enomoto works with right modules, whereas we are using left modules.

Fix some parameters $0 \not= \varepsilon_1 = q \in K$ and
$\varepsilon_2 = \varepsilon_3 = 1$.
Let $A_q = KQ/I_q$, where $Q$ is the Markov quiver
\[
\xymatrix@C=8mm@R=14mm{
& 2 \ar@<0.4ex>[dr]^{a_2}\ar@<-0.4ex>[dr]_{b_2}
\\
1 \ar@<0.4ex>[ur]^{a_1}\ar@<-0.4ex>[ur]_{b_1}
&& 3 \ar@<0.4ex>[ll]^{a_3}\ar@<-0.4ex>[ll]_{b_3}
}
\]
and $I_q$ is generated by the relations
\begin{align*}
\rho_{a_{i+2}} &:= 
a_{i+1}a_i - \varepsilon_i^{-1}b_{i+1}a_ib_{i+2}a_{i+1}b_i,
&b_{i+2}a_{i+1}a_i,
\\
\rho_{b_{i+2}} &:= b_{i+1}b_i - \varepsilon_i a_{i+1}b_ia_{i+2}b_{i+1}a_i,
&
a_{i+2}b_{i+1}b_i,
\end{align*}
with $i = 1,2,3$ and all indices interpreted modulo $3$.

\begin{Rem}
For $q =1$, $A_q$ is isomorphic to the Jacobian algebra associated to the once-punctured torus and the Labardini-Fragoso potential
$$
W_{\rm LF} := a_3a_2a_1 + b_3b_2b_1 - b_3a_2b_1a_3b_2a_1.
$$
\end{Rem}

Recall that an algebra $A$ is \emph{weakly symmetric} if
$P(i) \cong I(i)$ for all $i \in Q_0$.
In particular, such an $A$ is selfinjective.

\begin{Prop}[{\cite[Propositions~2.1 and 3.1]{E26}}]
\label{prop:eno1}
The algebra $A_q$ is a weakly symmetric $36$-dimensional
algebra with $\m^7 \subseteq I_q$, and we have
$$
\Omega_{A_q}^4(S) \cong \nu_{A_q}^{-1}(S) \cong S
$$
for all simple $A_q$-modules $S$.
\end{Prop}

(Enomoto does not explicitely say that $A_q$ is weakly symmetric, but it follows immediately from his calculations.)

\begin{Cor}
$A_q$ is $E$-symmetric.
\end{Cor}

\begin{proof}
Combine 
Proposition~\ref{prop:eno1} with
Corollary~\ref{cor:main2}.
\end{proof}

Thus by Theorem~\ref{thm:main1}, $A_q$ is a generalized Jacobian algebra.
Next, we construct a generalized hyperpotential for
$A_q$.

\begin{Prop}
Let $N_{I_q} := 7$.
Then  
$$
W := \{ \rho_{a_i}, \rho_{b_i} \mid i=1,2,3 \}
$$
is a generalized hyperpotential for $A_q$.
\end{Prop}

\begin{proof}
We have $\m^7 \subseteq I_q$.
We need to show that
$$
W^+ := W \cup Q_8
$$
generates the ideal $I_q$.

For $x,y \in KQ$, we write $x \equiv y$ if $x-y \in I(W)$.

We have
\begin{align*}
b_{i+2}a_{i+1}a_i &\equiv
\varepsilon_i^{-1} (b_{i+2}b_{i+1})a_ib_{i+2}a_{i+1}b_i
\\
&\equiv
\varepsilon_i^{-1}\varepsilon_{i+1}
(a_{i+2}b_{i+1}a_ib_{i+2}a_{i+1})a_ib_{i+2}a_{i+1}b_i.
\end{align*}
This implies 
$$
b_{i+2}a_{i+1}a_i \in I(W) + \m^9 \subseteq I(W) + \m^8.
$$
Similarly, we get $a_{i+2}b_{i+1}b_i 
\in I(W) + \m^8$.
We obtain
$$
I = I(W) + \m^8.
$$
In other words,
$W$ satisfies condition (h1) of Definition~\ref{def:genhyper}.

For a vertex $k$, let us vizualize the arrows incident to $k$:
$$
\xymatrix{
k+2 \ar@<0.4ex>[r]^<<<<{a_{k+2}} 
\ar@<-0.4ex>[r]_<<<<{b_{k+2}}
& k \ar@<0.4ex>[r]^<<<<{a_k}
\ar@<-0.4ex>[r]_<<<<{b_k} & k+1 
}
$$
(As before, all indices are interpreted modulo $3$.)
To display the maps $\scL_k$ and $\scR_k$ as matrices, we order the incoming and outgoing arrows as $(a_{k+2},b_{k+2})$ and 
$(a_k,b_k)$, respectively.
We have
$$
\scL_k,\; \scR_k\df P(k+2)^2 \to P(k+1)^2.
$$
For brevity, let
$$
p_k := a_{k+1}b_ka_{k+2}b_{k+1}
\text{\quad and \quad}
q_k := b_{k+1}a_kb_{k+2}a_{k+1}.
$$
Then
\begin{align*}
\scL_k &= 
\left(\begin{matrix}
\ov{l_{a_{k+2}}}(\rho_{a_k}) & \ov{l_{b_{k+2}}}(\rho_{a_k})
\\
\ov{l_{a_{k+2}}}(\rho_{b_k}) & \ov{l_{b_{k+2}}}(\rho_{b_k})
\end{matrix}\right)
=
\left(\begin{matrix}
a_{k+1} & -\varepsilon_{k+1}^{-1} p_k
\\
-\varepsilon_{k+1}q_k & b_{k+1} 
\end{matrix}\right),
\\
\scR_k &= 
\left(\begin{matrix}
\ov{r_{a_k}}(\rho_{a_{k+2}}) & \ov{r_{a_k}}(\rho_{b_{k+2}})
\\
\ov{r_{b_k}}(\rho_{a_{k+2}}) & \ov{r_{b_k}}(\rho_{b_{k+2}})
\end{matrix}\right)
=
\left(\begin{matrix}
a_{k+1} & -\varepsilon_k p_k
\\
-\varepsilon_k^{-1}q_k & b_{k+1} 
\end{matrix}\right),
\end{align*}
where the entries mean right multiplication by (the residue classes modulo $I$ of) the displayed elements.

Define
\begin{align*}
h_{k,{\rm in}} &:=
\left(\begin{matrix} \vep_k\vep_{k+1} & 0\\0&1 \end{matrix}\right)
\in \Aut_{A_q}(P(k+2)^2),
\\
h_{k,{\rm out}} &:= \left(\begin{matrix}
\vep_k\vep_{k+1}& 0\\0&1
\end{matrix}\right) \in \Aut_{A_q}(P(k+1)^2).
\end{align*}
We calculate
\begin{align*}
h_{k,{\rm out}}^{-1} \circ \scR_k \circ h_{k,{\rm in}} 
&=
\left(\begin{matrix}
\vep_k^{-1}\vep_{k+1}^{-1}& 0\\0&1
\end{matrix}\right)
\left(\begin{matrix}
a_{k+1} & -\varepsilon_k p_k
\\
-\varepsilon_k^{-1}q_k & b_{k+1} 
\end{matrix}\right)
\left(\begin{matrix}
\vep_k\vep_{k+1}& 0\\0&1
\end{matrix}\right)
\\
&= \left(\begin{matrix}
a_{k+1} & -\varepsilon_{k+1}^{-1} p_k
\\
-\varepsilon_{k+1}q_k & b_{k+1} 
\end{matrix}\right) = \scL_k.
\end{align*}
Thus $W$ also satisfies (h2).
\end{proof}

In the example above,
we have $\scL_k = \scR_k$ for each $k = 1,2,3$ if and only if $q=1$.

%%%%%%%%%%%%%%%%%%%%%%%%%%%%%%%%%%%%%%%%%%%
\subsection{Enomoto's example is not a Jacobian algebras}
%%%%%%%%%%%%%%%%%%%%%%%%%%%%%%%%%%%%%%%%%%%
For each $0 \not= \lambda \in K$, Enomoto studies the $2$-dimensional
indecomposable $A_q$-module
$U_\lambda$ corresponding to the following representation:
\[
\xymatrix@C=8mm@R=14mm{
& K \ar@<0.4ex>[dr]^{0}\ar@<-0.4ex>[dr]_{0}
\\
K \ar@<0.4ex>[ur]^{\lambda \cdot}\ar@<-0.4ex>[ur]_{1 \cdot}
&& 0 \ar@<0.4ex>[ll]^{0}\ar@<-0.4ex>[ll]_{0}
}
\]

\begin{Lem}[{\cite[Propositions~4.1]{E26}}]
\label{lem:eno2}
There is an exact sequence
$$
0 \to U_{q\lambda} \to P(2) \to P(1) \to U_\lambda \to 0.
$$
In particular, we have
$\Omega_{A_q}^{2m}(U_\lambda) \cong U_{q^m\lambda}$
for all $m \in \Z$.
\end{Lem}

\begin{Lem}\label{lem:eno3}
We have
$$
\nu_{A_q}(U_\lambda) \cong U_{\lambda/q}.
$$
\end{Lem}

\begin{proof}
This is not stated explicitely by Enomoto \cite{E26}, but it can be shown using his computations.
\end{proof}

\begin{Cor}\label{cor:eno5}
We have $\tau_{A_q}(U_\lambda) \cong U_\lambda$.
\end{Cor}

\begin{proof}
Since $A_q$ is selfinjective and $U_\lambda$ is indecomposable and non-projective, we have 
$$
\tau_{A_q}(U_\lambda) \cong \Omega_{A_q}^2(\nu_{A_q}(U_\lambda)).
$$
Now the claim follows from Lemmas~\ref{lem:eno2} and \ref{lem:eno3}.
\end{proof}

\begin{Cor}\label{cor:eno4}
For $q \not= 1$, $A_q$ is not stably
$3$-Calabi-Yau. 
\end{Cor}

\begin{proof}
The algebra $A_q$ is $1$-Iwanaga-Gorenstein
with $\GP(A_q) = \md(A_q)$, since
$A_q$ selfinjective.
We have 
$$
\Ext_{A_q}^1(U_\lambda,U_\lambda) \cong 
D\overline{\Hom}_{A_q}(U_\lambda,\tau_{A_q}(U_\lambda))
\cong D\overline{\End}_{A_q}(U_\lambda).
$$
Here we used the Auslander-Reiten formulas and Corollary~\ref{cor:eno5}.
Since $U_\lambda$ is non-injective and $\End_{A_q}(U_\lambda) \cong K$, we
get $\ext_{A_q}^1(U_\lambda,U_\lambda) = 1$.
On the other hand, we have
$$
\Ext_{A_q}^2(U_\lambda,U_{\lambda}) \cong
\underline{\Hom}_{A_q}(\Omega_{A_q}^2(U_\lambda),U_{\lambda}) 
\cong 
\underline{\Hom}_{A_q}(U_{q\lambda},U_{\lambda}).
$$
The first isomorphism holds, since $A_q$ is selfinjective, and the second isomorphism
comes from Lemma~\ref{lem:eno2}.
We get $\ext_{A_q}^2(U_\lambda,U_\lambda) = 0$, since
$\Hom_{A_q}(U_{q\lambda},U_\lambda) = 0$.
Thus $A_q$ is not stably $3$-Calabi-Yau.
\end{proof}

\begin{Cor}\label{cor:exnonjacobian}
For $q \not= 1$, $A_q$ is not isomorphic to a Jacobian algebra 
defined by a hyperpotential.
\end{Cor}

\begin{proof}
By Corollary~\ref{cor:eno4}, $A_q$ is not stably $3$-Calabi-Yau.
Now Theorem~\ref{thm:Jacobi3CY} implies that
$A_q$ cannot be a Jacobian algebra.
\end{proof}

\bigskip

%%%%%%%%%%%%%%%%%%%%%%%%%%%%%%%%%%%
\section*{Statement on AI usage}
%%%%%%%%%%%%%%%%%%%%%%%%%%%%%%%%%%%
This article relies strongly on the use of AI. 
It is based on \cite{BS26}, where we initiate the study of $E$-symmetric algebras and on \cite{GHLS26}, which gives 
general constructions of projective and injective presentations and relates these to $g$-vectors.
We asked the AI model ChatGPT-6 Astra for a generalization of the definition of a Jacobian algebra such that the resulting class of generalized Jacobian algebras coincides with the class of $E$-symmetric algebras.
We gave the findings from \cite{GHLS26} as an input, because we 
expected this to be the correct framework for our question. 
Amazingly, the AI came up with a suitable definition.
We took this definition and the suggested proof and transformed them
into an article
that can be understood by anyone who is familiar with basic Auslander-Reiten theory and the definition of Jacobian algebras (in the sense of Derksen--Weyman--Zelevinsky).
The AI also recognized that the AI-generated example from \cite{E26} is $E$-symmetric and not Jacobian.
All the writing was done by the three authors, who take full responsibility for the correctness of all results.

%%%%%%%%%%%%%%%%%%%%%%%%%%%%%%
\section*{Acknowledgements}
%%%%%%%%%%%%%%%%%%%%%%%%%%%%%%
CH and JS were funded by the Deutsche Forschungsgemeinschaft (DFG, German Research
Foundation) under Germany’s Excellence Strategy-GZ2047/1, Projekt-ID390685813.
JS thanks the
Faculty of Mathematics and Computer Science
of the Nicolaus Copernicus University in Toru\'n for 10 days of hospitality in September 2026.

%%%%%%%%%%%%%%%%%%%%%%%%%%%%%%%%%%%%%%%%%%%%%

\end{document}